# SUBSTRATA IN REDUCTIVE GROUPS

G. Lusztig

**1.** Let $G$ be a connected reductive algebraic group over an algebraic closed field $\mathbf{k}$ of characteristic $p \geq 0$. Let $W$ be the Weyl group of $G$; let $Rep(W)$ be the set of finite dimensional representations of $W$ over $\mathbf{Q}$ (up to isomorphism) and let $\mathrm{Irr}(W)$ be the subset of $Rep(W)$ consisting of irreducible representations. Let $\delta : G \to 2\mathbf{N}$ be the function which to any $g \in G$ associated the dimension of the conjugacy class of $g$ in $G$. (It is known [Sp82] that this dimension is always even.) For any $d \in 2\mathbf{N}$ let $G^d = \delta^{-1}(d)$, that is the union of all conjugacy classes of $G$ that have dimension $d$. According to [L15,0.1],

(a) *the set $\{d \in 2\mathbf{N}; G^d \neq \emptyset\}$ depends only on $W$ (and not on $\mathbf{k}$).*

In [L15] a partition of $G$ into finitely many subsets (called strata) was defined. Each stratum is a union of conjugacy classes contained in a fixed $G^d$. The strata of $G$ are parametrized by a subset $\mathrm{Irr}_*(W)$ of $\mathrm{Irr}(W)$ which depends only on the type of $G$ and not on $\mathbf{k}$.

**2.** In this paper we define a partition of $G$ into finitely many subsets (called substrata). Each substratum is a union of conjugacy classes contained in a fixed stratum of $G$. The substrata of $G$ are parametrized by a subset $Rep_*(W)$ of $Rep(W)$. We will also formulate a number of conjectural properties of substrata.

**3.** Let $\mathcal{B}$ be the variety of Borel subgroups of $G$. Let $\nu = \dim \mathcal{B}$. For any $g \in G$ let $\mathcal{B}_g = \{B \in \mathcal{B}; g \in B\}$ be the Springer fibre at $g$. It is known [Sp82] that $\mathcal{B}_g$ is nonempty of pure dimension $b(g) = \nu - \delta(g)/2$. Let $l$ be a prime number $\neq p$. In [L81] an action of $W$ on the $l$-adic homology $H_i(\mathcal{B}_g, \bar{\mathbf{Q}}_l)$ was defined for any $i$, using intersection cohomology. (In the case where $g$ is unipotent and $p \gg 0$ or $p = 0$, this action was first defined by Springer [Sp76] using a different method.) This action comes from a well defined ${}^i V_g \in Rep(W)$. We write $V_g = {}^{2b(g)} V_g$.

We define a $W$-submodule $V'_g$ of $V_g$ by descending induction on $b(g)$. We have $b(g) \leq \nu$. Assume first that $b(g) = \nu$. Then $V'_g = V_g$ (this is the sign representation of $W$). Next we assume that $b(g) < \nu$ and that $V'_{g'}$ is already defined when $g' \in G, b(g') > b(g)$. Let $V'_g$ be the sum of all irreducible $W$-submodules $E$ of $V_g$ such that $E$ does not appear in $V'_{g'}$ for any $g' \in G$ such that $b(g') > b(g)$. This completes the inductive definition of $V'_g$.

Given $g_1, g_2$ in $G$ we say that $g_1 \sim g_2$ if $V'_{g_1} = V'_{g_2}$. This is an equivalence relation on $G$; the equivalence classes are said to be the substrata of $G$.

For $g \in G$, the image of the natural map $H_{2b(g)}(\mathcal{B}_g, \bar{\mathbf{Q}}_l) \to H_{2b(g)}(\mathcal{B}, \bar{\mathbf{Q}}_l)$ can be regarded as an object $V''_g$ of $Rep(W)$; it is actually in $\text{Irr}(W)$ and even in $\text{Irr}_*(W)$.

For any $E \in \text{Irr}(W)$ let $n_E$ be the smallest integer $i$ such that $E$ appears in the $W$-module $H_{2i}(\mathcal{B}, \bar{\mathbf{Q}}_l)$. It is known (see [L15, 2.2]) that

(a) $V''_g$ appears with multiplicity one in $V_g$ and $n_{V''_g} = b(g)$; moreover all other composition factors $E'$ of $V_g$ satisfy $n_{E'} > b(g)$.

We show:

(b) $V''_g$ appears with multiplicity one in $V'_g$ and all other composition factors $E'$ of $V'_g$ satisfy $n_{E'} > b(g)$.

By (a), it is enough to show that $V''_g$ appears in $V'_g$ or equivalently that $V''_g$ does not appear in $V'_{g'}$ if $b(g') > b(g)$. If it did appear then it would also appear in $V_{g'}$. By (a), this would imply $n_{V''_g} \geq b(g')$ hence $n_{V''_g} > b(g')$, so that $b(g) > b(g')$, a contradiction.

(c) $V''_g$ is uniquely determined by $V'_g$. Moreover, if $V''_{g'}$ appears in $V'_g$ then $V''_{g'} = V''_g$.

The first statement follows from (b). Let $E = V''_{g'}$. If $E$ appears in $V'_g$ and $E \neq V''_g$ then by (b) we have $n_E > b(g)$. We have also $n_E = b(g')$ hence $b(g') > b(g)$ and $E$ appears in $V'_{g'}$. This contradicts the definition of $V'_g$. Thus (c) holds.

Recall that $g_1, g_2$ in $G$ are said to be in the same stratum if $V''_{g_1} = V''_{g_2}$. Using (c) we deduce the following result.

(d) If $g_1, g_2$ in $G$ are in the same substratum then they are in the same stratum.

Let $Rep_*(W)$ be the subset of $Rep(W)$ formed by the $W$-modules $V'_g$ for various $g \in G$. Note that

(e) $Rep_*(W)$ is a finite set.

Indeed from [L15,2.2] we see that the $W$-modules of the form $V_g$ ($g \in G$) form a finite subset of $Rep(W)$. Since there are only finitely many objects of $Rep(W)$ which are isomorphic to $W$-submodules of $V_g$ we see that (e) holds.

Note that $Rep_*(W)$ can be viewed as an indexing set of the set of substrata of $G$. Hence the number of substrata of $G$ is finite.

We have a partition $Rep_*(W) = \sqcup_\Sigma Rep_\Sigma(W)$ where $\Sigma$ runs over the strata of $G$ and $Rep_\Sigma(W)$ is the set of all objects in $Rep_*(W)$ that contain $V''_g$ where $g \in \Sigma$. (This is a partition by (c).)

**4.** We define a (surjective) map

(a) $\xi : \text{Irr}(W) \to \text{Irr}_*(W)$

by $E \mapsto \underline{\tau}_0(\emptyset, E, \bar{\mathbf{Q}}_l)$ (notation of [L23, 1.8]).

We now state a number of conjectural properties of substrata.

(b) *Let $\Sigma$ be a stratum of $G$ and let $E_\Sigma$ be the corresponding object of $\text{Irr}_*(W)$. If $E \in \text{Irr}(W)$ appears in a $W$-module in $Rep_\Sigma(W)$ then $\xi(E) = E_\Sigma$. Moreover, the $W$-modules in $Rep_\Sigma(W)$ are related to the $W$-modules in $\xi^{-1}(E_\Sigma)$ by a triangular* matrix with diagonal entries 1. In particular $Rep_\Sigma(W)$ is in canonical bijection with $\xi^{-1}(E_\Sigma)$. Hence $Rep_*(W)$ is in canonical bijection with $\mathrm{Irr}(W)$.

(c) For any stratum $\Sigma$, $Rep_\Sigma(W)$ depends only on $W$ and not on $\mathbf{k}$. Hence $Rep_*(W)$ depends only on $W$ and not on $\mathbf{k}$.

**5.** We now assume that $G = GL_n(\mathbf{k})$. For any partition $n_* = (n_1 \leq n_2 \leq \cdots \leq n_k)$ of $n$ let $X_{n_*}$ be the set of semisimple elements of $G$ whose centralizer is isomorphic to $GL_{n_1}(\mathbf{k}) \times GL_{n_2}(\mathbf{k}) \times \ldots \times GL_{n_k}(\mathbf{k})$. We have $X_{n_*} \subset G^d$ where $d = n^2 - (n_1^2 + n_2^2 + \cdots + n_k^2)$. Let $\bar{X}_{n_*}$ be the closure of $X_{n^*}$ in $G^d$ (see no.1). It is known (see [C20]) that the subsets $\bar{X}_{n_*}$ are exactly the strata of $G$. They are also the substrata of $G$. We have $Rep_*(W) = \mathrm{Irr}_*(W) = \mathrm{Irr}(W)$. We see that in this case 4(b),(c) hold.

**6.** We now assume that $G = Sp_4(\mathbf{k})$. In this case we have $Rep_*(W) = \mathrm{Irr}_*(W) = \mathrm{Irr}(W)$; the substrata of $G$ are the same as the strata of $G$. They are the 5 subsets $G^8, G^6, G_I^4, G_{II}^4, G^0$ where $G^d$ are as in no.1 and

$G_I^4$ consists of the conjugacy class of a transvection $T$ union with (if $p \neq 2$) that of $-T$;

$G_{II}^4$ consists of the conjugacy class of a unipotent element $\tilde{T}$ in $G^4$ not conjugate to $T$ (if $p = 2$) and of the unique semisimple class in $G^4$ (if $p \neq 2$).

We have $G^4 = G_I^4 \sqcup G_{II}^4$.

We see that in this case 4(b),(c) hold.

**7.** We now assume that $G = Sp_6(\mathbf{k})$. Now $\mathrm{Irr}(W)$ consists of 10 objects denoted

$$1_0, 3_1, 2_2, 3_2, 1_3, 3_3, 3_4, 2_5, 1_6, 1_9.$$

Here $d_n$ denotes a representation $E \in \mathrm{Irr}(W)$ of dimension $n$ and such that $n_E = n$.

The following subsets of $G^{14}$ (with $G^d$ as in no.1) form a partition of $G^{14}$.

(a) The conjugacy class of a unipotent element with Jordan blocks $1 + 1 + 4$, union (if $p \neq 2$) with $-1$ times this conjugacy class.

(b) The conjugacy class of a unipotent element with Jordan blocks $3 + 3$, union (if $p \neq 2$) with $-1$ times this conjugacy class.

(c) The union of conjugacy classes of semisimple elements with centralizer $GL_2(\mathbf{k}) \times Sp_2(\mathbf{k})$.

(d) (if $p \neq 2$) The union of conjugacy classes of semisimple elements with centralizer $Sp_2(\mathbf{k}) \times Sp_2(\mathbf{k}) \times \mathbf{k}^*$.

(e) The union of conjugacy classes of elements $su$ where $s \in G$ is semisimple with centralizer $Sp_4(\mathbf{k}) \times \mathbf{k}^*$ and $u$ is $T \times 1$ ($T$ as in no.6).

(f) (if $p = 2$) The union of conjugacy classes of elements $su$ where $s \in G$ is semisimple with centralizer $Sp_4(\mathbf{k}) \times \mathbf{k}^*$ and $u$ is $\tilde{T} \times 1$ ($\tilde{T}$ as in no.6).

(g) (if $p \neq 2$) The union of conjugacy classes of elements $su$ where $s \in G$ is semisimple with centralizer $Sp_4(\mathbf{k}) \times Sp_2(\mathbf{k})$ and $u = u' \times 1$ where $u' \in Sp_4(\mathbf{k})$ is subregular unipotent.

(h) (if $p \neq 2$) The union of conjugacy classes of elements $su$ where $s \in G$ is semisimple with centralizer $Sp_4(\mathbf{k}) \times Sp_2(\mathbf{k})$ and $u$ is $T \times u'$ ($T$ as in no.2, $u'$ regular unipotent).

Let $G_I^{14}$ be the union of the sets (a),(e) (if $p = 2$) or (a),(e),(h) (if $p \neq 2$). Let $G_{II}^{14}$ be the set (c) (if $p = 2$) or the union of the sets (b),(c) (if $p \neq 2$). Let $G_{III}^{14}$ be the union of the sets (b),(f) (if $p = 2$), or (d),(g),(h) (if $p \neq 2$).

The set $Rep_*(W)$ consists of 10 objects

$$1_0, 3_1, 2_2, 3_2, 3_2 + 1_3, 3_3, 3_4, 2_5, 1_6, 1_9.$$

The corresponding 10 substrata of $G$ are

$$G^{18}, G^{16}, G_I^{14}, G_{II}^{14}, G_{III}^{14}, G^{12}, G^{10}, G^8, G^6, G^0.$$

The set $\mathrm{Irr}_*(W)$ consists of 9 objects:

$$1_0, 3_1, 2_2, 3_2, 3_3, 3_4, 2_5, 1_6, 1_9.$$

The corresponding 9 strata of $G$ are

$$G^{18}, G^{16}, G_I^{14}, G_{II}^{14} \sqcup G_{III}^{14}, G^{12}, G^{10}, G^8, G^6, G^0.$$

We see that in this case 4(b),(c) hold.

## References

[C20] G.Carnovale, *Lusztig strata are locally closed*, Arch.Math. (Basel) **115** (2020), 23-26.
[L81] G.Lusztig, *Green polynomials and singularities of unipotent classes*, Adv.in Math. **42** (1981), 169-178.
[L15] G.Lusztig, *On conjugacy classes in a reductive group*, Representations of Reductive Groups, Progr.in Math. 312, Birkhäuser, 2015, pp. 333-363.
[L23] G.Lusztig, *Unipotent character sheaves and strata of a reductive group*, Repres.Th. **27** (2023), 1126-1141.
[Sp82] N.Spaltenstein, *Classes unipotentes et sous-groupes de Borel*, Lecture Notes in Math. 946, Springer Verlag, 1982.
[S76] T.A.Springer, *Trigonometric sums, Green functions and representations of Weyl groups*, Invent.Math. **36** (1976), 173-207.

Department of Mathematics, M.I.T., Cambridge, MA 02139